\documentclass{article}

\usepackage{arxiv}

\usepackage[utf8]{inputenc} 
\usepackage[T1]{fontenc}    
\usepackage{hyperref}       
\usepackage{url}            
\usepackage{booktabs}       
\usepackage{amsfonts}       
\usepackage{nicefrac}       
\usepackage{microtype}      
\usepackage{graphicx}
\usepackage{natbib}
\usepackage{doi}
\usepackage{bm}
\usepackage{makecell}
\usepackage{multirow}
\usepackage[textsize=tiny]{todonotes}
\usepackage{xcolor, soul}
\usepackage{threeparttable}
\usepackage{adjustbox}
\usepackage{siunitx} 
\usepackage{textcomp}
\usepackage{gensymb}
\usepackage{tikz}

\usepackage[c1]{optidef}

\title{Direct multiple shooting and direct collocation perform similarly in biomechanical predictive simulations}

\author{ \href{https://orcid.org/0000-0002-9335-630X}{\includegraphics[scale=0.06]{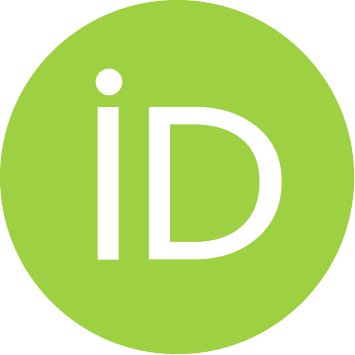}\hspace{1mm}P. Puchaud} \\
	Laboratoire de Simulation et Modélisation du Mouvement \\
    École de kinésiologie et des Sciences de l’Activité Physique \\
    Université de Montréal, Laval, Canada\\
    Centre de réadaptation Marie-Enfant \\
    CHU Sainte-Justine, Montréal, Canada \\
	\texttt{pierre.puchaud@umontreal.ca} \\
	\And
	\href{https://orcid.org/0000-0000-0000-0000}{\includegraphics[scale=0.06]{orcid.pdf}\hspace{1mm}F. Bailly} \\
	INRIA, Univ Montpellier, Montpellier, France \\
	\texttt{francois.bailly@inria.fr} \\
	\AND
	\href{https://orcid.org/0000-0002-9335-630X}{\includegraphics[scale=0.06]{orcid.pdf}\hspace{1mm}M. Begon} \\
	Laboratoire de Simulation et Modélisation du Mouvement \\
    École de kinésiologie et des Sciences de l’Activité Physique \\
    Université de Montréal, Laval, Canada\\
    Centre de réadaptation Marie-Enfant \\
    CHU Sainte-Justine, Montréal, Canada \\
	\texttt{mickael.begon@umontreal.ca} \\
}

\renewcommand{\shorttitle}{DMS and DC perform similarly in biomechanical predictive simulations}

\DeclareFontFamily{U}{mathb}{\hyphenchar\font45}
\DeclareFontShape{U}{mathb}{m}{n}{
      <5> <6> <7> <8> <9> <10> gen * mathb
      <10.95> mathb10 <12> <14.4> <17.28> <20.74> <24.88> mathb12
      }{}
\DeclareSymbolFont{mathb}{U}{mathb}{m}{n}
\DeclareFontSubstitution{U}{mathb}{m}{n}

\let\dot\relax
\DeclareMathAccent{\dot}{0}{mathb}{"39}
\let\ddot\relax
\DeclareMathAccent{\ddot}{0}{mathb}{"3A}
\let\dddot\relax
\DeclareMathAccent{\dddot}{0}{mathb}{"3B}
\let\ddddot\relax
\DeclareMathAccent{\ddddot}{0}{mathb}{"3C}

\newcommand{\qddot}{\bm{\ddot{q}}}
\newcommand{\qdot}{\bm{\dot{q}}}
\newcommand{\q}{\bm{q}}
\newcommand{\am}{\bm{\sigma}}
\newcommand{\btau}{\bm{\tau}}
\newcommand{\controls}{\bm{\mathrm{u}}}
\newcommand{\states}{\bm{\mathrm{x}}}
\newcommand{\sysd}{\bm{\mathrm{f}}}
\newcommand{\sinit}{\bm{\mathrm{s}}_0}
\newcommand{\statesdot}{\bm{\dot{\mathrm{x}}}}

\newcommand{\ERK}{\mbox{\textit{MSE}}}
\newcommand{\IRKid}{\mbox{\textit{MSI}\textsuperscript{ID}}}
\newcommand{\IRKfd}{\mbox{\textit{MSI}\textsuperscript{FD}}}
\newcommand{\DCfd}{\mbox{\textit{DC}\textsuperscript{FD}}}
\newcommand{\DCid}{\mbox{\textit{DC}\textsuperscript{ID}}}

\newcommand{\knotpoints}{\mbox{knot points}}

\newcommand{\htmlversion}{A Html version of the figure is available as supplementary material to further explore the data.}

\hypersetup{
pdftitle={Direct multiple shooting and direct collocation show similar performances in biomechanical predictive simulations},
pdfauthor={Pierre Puchaud, François Bailly, Mickael Begon},
}

\begin{document}
\maketitle

\begin{abstract}
Direct multiple shooting (DMS) and direct collocation (DC) are two common transcription methods for solving optimal control problems (OCP) in biomechanics and robotics.  
They have rarely been compared in terms of solution and speed. Through five examples of predictive simulations solved using five transcription methods and 100 initial guesses in the Bioptim software,  we showed that not a single method outperformed systematically better.
All methods converged to almost the same solution (cost, states, and controls) in all but one OCP, with several local minima being found in the latter.  
Nevertheless, DC based on fourth-order Legendre polynomials provided overall better results, especially in terms of dynamic consistency compared to DMS based on a fourth-order Runge-Kutta method. Furthermore, expressing the rigid-body constraints using inverse dynamics was usually faster than forward dynamics.  
DC with dynamics constraints based on inverse dynamics converged to better and less variable solutions.  
Consequently, we recommend starting with this transcription to solve OCPs but keep testing other methods. 
\end{abstract}


\section{Introduction}

Predictive simulation based on the optimal control theory can be insightful in biomechanics and robotics. In both fields, rigid-multibody dynamics equations describe the physical behavior of the system. In robotics, optimal control is used in trajectory planning \citep{Diehl2006FastControl}.
In biomechanics, human motor behavior often incorporates elements of optimality, and optimal solutions can give us insight into how humans plan and carry out movements, as well as provide innovative techniques for solving complex tasks.
It is worth noting that the models in biomechanics can rely on joint torque-driven models, as for robotics, or muscle-driven models. The latter adds an extra layer of dynamics to optimal control programs (OCPs), leading typical musculoskeletal models to handle thousands of variables, and enlarging the OCPs \citep{Febrer-Nafria2022PredictiveReview}. 
Direct approaches are commonly used for transcribing OCPs related to rigid-body dynamics, as opposed to indirect approaches, which tend to be more challenging to implement and initialize \citep{Kelly2017AnCollocation, Betts2010PracticalProgramming}. 
The most common transcription methods are direct single shooting, direct multiple shooting (DMS), and direct (DC) collocation. 
The superiority of one transcription over another is still a topic of ongoing debate. While DMS is more popular in robotics \citep{Giftthaler2018AControl, Sleiman2019Contact-ImplicitManipulation}, DC has gained popularity in biomechanics \citep{Ezati2019AGait, Febrer-Nafria2020ComparisonModel, Lee2016GeneratingMATLAB}. It is important to compare these transcriptions for the biomechanics community.

Direct single shooting has been used in biomechanics until recently \citep{Geijtenbeek2019SCONE:Motion}. It has limitations, such as being susceptible to local minima and chaotic behavior due to its dependency on the initial solution. The direct single shooting approach inherits these ill-conditioned OCPs that are often not solvable numerically speaking and slow \citep{Diehl2006FastControl}.
DMS and DC approaches improve numerical stability and convergence time compared to direct single shooting \citep{Betts2010PracticalProgramming}.
DMS uses time-stepping numerical integration to simulate consecutive intervals of a trajectory, and the arrival node of intervals is constrained to be equal to the first node of the next interval \citep{Bock1984AProblems}.  
DMS transcription has proven its efficiency in solving various biomechanical problems, including tracking motions \citep{Venne2022}, predictive simulations such as acrobatic moves \citep{Koschorreck2012ModelingTwists, Charbonneau2020OptimalTrampoline}, gait \citep{Felis2016SynthesisMethods}, lifting techniques with an exoskeleton \citep{Sreenivasa2017OptimalModel} and its optimal settings \citep{Harant2017ParameterControl, Sreenivasa2018PredictingControl}.
In DC, the system states are approximated with polynomial splines on a temporal grid, and optimization ensures that the splines obey the system's dynamics at collocation points \citep{Biegler1984SolutionCollocation}. 
Direct collocation has been used for different biomechanics-related problems such as tracking motions \citep{Lin2017Three-dimensionalCollocation}, and gait predictive simulations \citep{Ackermann2010OptimalityGait, Miller2015OptimalRunning, Falisse2022ModelingWalking}.
However, both DMS and DC approaches have been applied to solve biomechanical problems, but no study has provided a fair comparison of these transcriptions.

DMS and DC transcriptions are not that different when looking at their transcriptions (see Sec.~\ref{back:ocp}). 
Indeed, when using a DMS approach with an implicit Runge-Kutta (IRK) solver, the splines obey the system's dynamics at collocation points through a Newton descent in the IRK solver. 
When using direct collocation, the same splines obey the dynamics through the constraints of the non-linear program (NLP) solver. 
The implicit formulations used to build the rigid-body dynamic constraints ($F(\states, \statesdot, \controls) = 0$) affects the results \citep{Ferrolho2021InverseOptimization}.
Researchers have reported key methodological advances, including efficiently handling rigid-body and muscle dynamics through implicit formulations \citep{VanDenBogert2011ImplicitControl, DeGroote2016, Puchaud2023OptimalityDynamics}, which have been recommended to solve OCPs in biomechanics. 
For instance, \cite{Docquier2019ComparisonProblem} showed that direct collocations OCPs reached lower optima when rigid-body dynamic constraints matched joint acceleration through forward dynamics ($\qddot - FD(\q, \qdot, \btau) = 0$) and that the CPU time was lower matching joint torques through inverse dynamics ($\btau - ID(\q, \qdot, \qddot) = 0$).
In an effort to include as many formulations as possible, the present study will assess the performance between these transcriptions (DMS with explicit RK, DMS with IRK, DC) and the formulation of rigid-body dynamic constraints.

Evaluating the optimality of OCPs commonly involves sensitivity studies with initial guesses.
OCP have either been started with data-informed initial guesses close to an expected solution \citep{Falisse2019} or with noise \citep{Charbonneau2022OptimalSelf-collision}.
For the present study, multi-start approaches with random noise were chosen because they are more challenging as they test the ability of OCPs to converge toward the same optimum \citep{Puchaud2023OptimalityDynamics}. Such approaches are needed to compare different transcriptions between them in terms of optimality.
Recently, more and more direct transcriptions rely on fixed-step ordinary differential equations (ODE) integrators, favored by algorithmic differentiation \citep{Andersson2019CasADi:Control}. Unfortunately, the integration errors cannot be controlled anymore by the integrators while solving an OCP by using adaptive step size ODE solvers \citep{Serban2008CVODES:SUNDIALS} or adaptive
mesh refinement \citep{Patterson2014GPOPS-II}. Thus, there is a need for metrics to assess the dynamic consistency of the optimal solutions found. While this step is essential to evaluate the dynamics consistency of the solution, this practice has not been frequently done when using direct collocation \citep{Michaud2022BioptimBiomechanics}. The dynamics consistency of direct transcriptions must be examined to estimate the reliability of all optimal solutions found with multi-start approaches.

This study aims to evaluate the performance of various direct optimal control transcriptions and dynamic formulations on five OCPs.
First, we will review how the optimal control problem is transcribed from the continuous optimal control problem into non-linear programming through DC and DMS transcriptions.
Second, we will present two ways of formulating rigid-body constraints through forward or inverse dynamics.
Then, five problems of increasing complexity will be presented to benchmark these transcriptions. 
To evaluate the performance of each transcription, the convergence rates, the CPU time to converge, optimality, and dynamic consistency of solutions will be evaluated in a multi-start scheme.

\section{Background} \label{sec:Background}

This section will present the transcription of continuous trajectory optimization problems (\ref{back:ocp}) into direct collocation (\ref{back:dc}) and direct multiple shooting (\ref{back:dms}). Finally, the dynamics constraints specific to rigid-body dynamics will be described (\ref{back:dynamicconstraint}).

\subsection{Trajectory optimization} \label{back:ocp}

Let's consider any system whose dynamics is governed by an ODE of the form $\statesdot = \sysd(\states, \controls)$,
with $\states$ the state vector of size $n_{\states}$ and $\controls$ the control vector of size $n_{\controls}$. 
A generic trajectory optimization problem can be written as:

\small
\begin{subequations}
\label{eq:ocpcontinuous}
\begin{eqnarray}
\hspace{-5em} 
\underset{\substack{\states(\cdot),\controls(\cdot)}}{\min }
\ \ \
    & & \hspace{-2.5em}   \phi(\states(t),\controls(t)) \label{eq:cost_continuous} \\
    \text{s.t.}& & \states(0) - \sinit = 0,  \quad  \textit{Initial state constraint} \label{eq:initial_constraint_continuous} \\
    & \forall \; t \in [0,T], &\dot{\states}(t) - \sysd(\states(t),\controls(t)) = 0, \quad \textit{Dynamics constraints} \label{eq:dynamic_constraint_continuous} \\
    &  \forall \; t \in [0,T], &\bm{g}(\states(t),\controls(t)) \leq 0, \quad \textit{Path constraints} \label{eq:control_path_constraint_continuous}\\
    & &\bm{r}(\states(T)) \leq 0, \quad \textit{Terminal constraint} \label{eq:terminal_contraints_continuous}
\end{eqnarray}
\end{subequations}
\normalsize

\noindent with $\bm{s}_0$ the initial state and $T$ the final time of the problem. 
The cost function is denoted $\phi$. 
The inequality constraints are gathered in a function $\bm{g}$, and the terminal constraint is represented by $\bm{r}$.
This continuous-time optimization problem must be transcribed into a finite-dimensional nonlinear program (NLP) to be numerically tractable.
To discretize this problem, two transcription approaches are described in the following: \textit{direct collocation}, and \textit{direct multiple shooting}. The main differences between them are the way they handle the simulation and optimization. DC does both simultaneously, whereas DMS first simulates and then, optimizes the trajectories.

\subsection{Direct collocation transcription} \label{back:dc}

In direct collocation, the control trajectory $\controls(t)$ is discretized into a finite-dimension vector $\controls=[\controls_0 \ldots \controls_n \ldots \controls_{N-1}]$ of size $N$, at times~$t_n$, called \knotpoints{} with $n \in [0..N-1]$.
The state trajectory $\states(t)$ is similarly discretized with
$\states=[\states_0 \ldots \states_n \ldots \states_N]$ of size $N + 1$, at~times~$t_n$, with $n \in [0..N]$.  
$M$ collocation points subdivide the intervals $[t_n, t_{n+1}]$ at times $t_{n,m}$, with $m\in [0..M-1]$.
At these collocation points, intermediate states are denoted by $\states_{n,m}$.
For the sake of compactness, at each \knotpoints{} we introduce a vector of intermediate states $\Tilde{\states}_n = [\states_{n,0}, \ldots, \states_{n,M-1}]$.
The main idea behind collocation transcriptions is to model the state trajectory $\states(t)$ by a set of $N$ polynomials $\bm{p}_n$, parametrized by the intermediate states $\states_{n,m}$.
In the present work, we will focus on the specific case of \textit{orthogonal collocations}, where the polynomials are Lagrange ones, such that:
\small
\begin{equation} \label{eq:lagrange}
\begin{aligned}
    \forall \; t \in [t_n, t_{n+1}], \quad \bm{p}_n(t,\Tilde{\states}_n) = \sum_{m=0}^{M} \states_{n,m} \; \ell_{n,m}(t) \; \in \mathbb{R}^{n_{\states}}, \\
    \ell_{n,m}(t) = \prod_{k=0,p\ne m}^M \frac{t-t_{n,k}}{t_{n,m}-t_{n,k}} \in \mathbb{R}.
\end{aligned}
\end{equation}
\normalsize
Properties of Lagrange polynomials entail that $\bm{p}_n(t_{n,m},\Tilde{\states}_n) = \states_{n,m}$. Then, the dynamics constraints (Eq.~\ref{eq:dynamic_constraint_continuous}) are satisfied at every collocation point:

\small
\begin{equation} \label{eq:Colloc-Conditions}
\begin{aligned}
\forall \; (n,m) \in [0..N-1]\!\times\![0..M-1], \quad \sysd\left(\bm{p}_n(t_{n,m}, \Tilde{\states}_{n})  ,\; \controls_n\right) &= \dot{\bm{p}}_n(t_{n,m}, \Tilde{\states}_{n}).
\end{aligned}
\end{equation}
\normalsize

These collocation conditions are included in the NLP, yielding a large-scale but sparse problem, which can be written as:

\small
\begin{subequations}
\label{eq:ocpDC}
\begin{eqnarray}
\hspace{-5em} 
\underset{\substack{\Tilde{\states},\controls}}{\min }
\ \ \
    & \hspace{-2.5em}   \phi(\Tilde{\states},\controls) \label{eq:cost_DC} \\
    \text{s.t.}& &   \bm{\states}_{0,0} - \bm{s}_0 = 0,  
    \label{eq:initial_constraint_DC} \\
    & \forall \; (n, m)  \in [0..N-1]\!\times\![0..M-1], \; & \sysd\left(\bm{p}_n(t_{n,m}, \Tilde{\states}_{n})  ,\; \controls_n\right) = \dot{\bm{p}}_n(t_{n,m}, \Tilde{\states}_{n}),
    \label{eq:collocation_condition_DC} \\
    &\forall \; n \in [0..N-1], & \bm{p}_n(t_{n+1},\Tilde{\states}_n) - \states_{n+1,0} = 0,
    \label{eq:continuity_conditions_DC} \\
    &  \forall \; n \in [0..N-1],  & \bm{g}(\states_n,\controls_n) \leq 0,
    \label{eq:inequality_contraints_DC} \\
    &  & \bm{r}(\states_N) \leq 0.
    \label{eq:terminal_contraints_DC}
\end{eqnarray}
\end{subequations}
\normalsize
One may note that the optimization variable $\Tilde{\states}$ directly corresponds to the coefficients of the Lagrange polynomials.
Continuity constraints on $\bm{p}_n$ (Eq.~\ref{eq:continuity_conditions_DC}) ensure state continuity between the end of one interval and the beginning of the next one. 

\subsection{Direct multiple shooting transcription} \label{back:dms}
We introduce the discretized control trajectory
$\controls=[\controls_0 \ldots \controls_n \ldots \controls_{N-1}]$ of size $N$ and the discretized state trajectory $\states=[\states_0 \ldots \states_n \ldots \states_N]$ of size $N + 1$.
In DMS, the system ODE (Eq.~\ref{eq:dynamic_constraint_continuous}), is numerically solved on each interval as an initial value problem. State continuity constraints entail that the end of one interval coincides with the beginning of the next, to ensure the continuity of the trajectory, leading to the following NLP:

\small
\begin{subequations}
\label{eq:dms}
\begin{eqnarray}
\hspace{-5em} 
\underset{\substack{\states,\controls}}{\min }
\ \ \
    & \hspace{-2.5em}   \phi(\states,\controls) \label{eq:cost_dms} \\
    \text{s.t.}& &  \bm{\states}_{0,0} - \bm{s}_0 = 0,  
    \label{eq:dms_initial_constraint} \\
    &  \forall \; n  \in [0..N-1], & \states_{n+1} = \states_{n} + \int_{t_n}^{t_{n+1}} \sysd(\states_n,\controls_n) \; dt,
    \label{eq:dms_continuity_conditions} \\
    &  \forall \; n \in [0..N-1],  & \bm{g}(\states_n,\controls_n) \leq 0,
    \label{eq:dms_inequality_contraints} \\
    &  & \bm{r}(\states_N) \leq 0, 
    \label{eq:dms_terminal_contraints}
\end{eqnarray}
\end{subequations}
\normalsize

The flexibility of DMS comes from the choice of the ODE solver needed to compute the continuity constraint (Eq.~\ref{eq:dms_continuity_conditions}), which ensures the consistency of the dynamics, for example: 
\begin{itemize}
    \item explicit Runge-Kutta (ERK) integrators (\textit{e.g.}, 2nd order, 4th order), which are fast and precise enough for most applications;
    \item collocation-based IRK; interestingly, the collocation conditions of the DC transcription (Eq. \ref{eq:Colloc-Conditions}) are mathematically equivalent to a DMS formulation with an implicit Runge-Kutta ODE solver \citep{Wright1970SomeProperties, Serie1969LaCollocation, Hairer2006NumericalIntegrators}. Numerically, the DMS transcription with IRK ensures the constraints through a Newton-Raphson descent, while, in DC, they are solved directly in the optimizer;
    \item variable-step and order-step integrators, such as RK45 or CVODES \citep{Serban2008CVODES:SUNDIALS}, make this transcription method theoretically as precise as desired regarding the dynamical consistency of the solutions, thanks to the adaptive step size control. 
\end{itemize}

\begin{figure*}[h!]
\includegraphics[width=1\linewidth,trim={0 0 0 0}, clip]{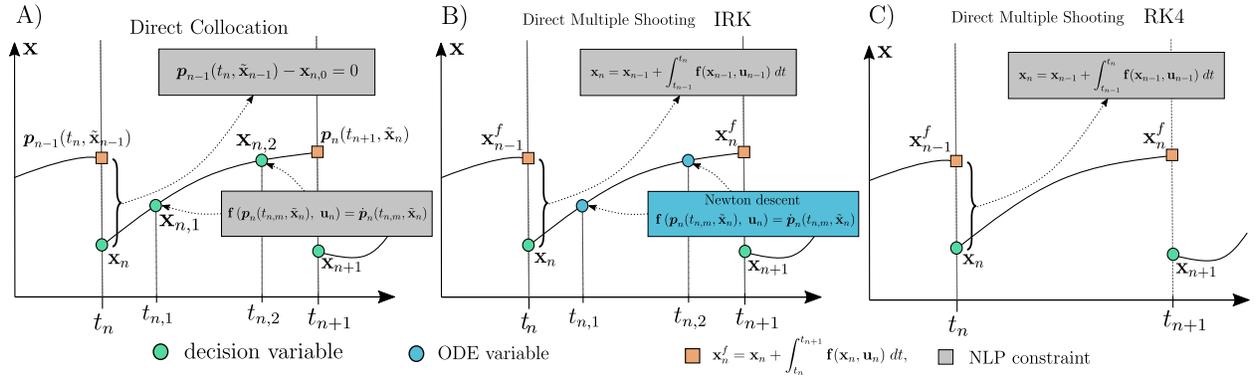}
\centering
\caption{Illustrations of OCP transcriptions. A) Direct Collocation, B) Direct multiple shooting with IRK, C) Direct multiple Shooting. }
\label{fig:transcriptions}
\end{figure*}

\subsection{Rigid-body dynamics constraints} \label{back:dynamicconstraint}

In the previous sections, we have introduced two transcriptions methods for trajectory optimization, namely, DC and DMS. 
As presented, they can be declined into several formulations (changing polynomials for DC, changing integrators for DMS), and some parallels can be drawn between them.
In the following, we propose to investigate the behavior of a subset of selected formulations, namely: (\textit{i}) DMS with explicit ERK, (\textit{ii}) DMS with IRK, and (\textit{iii}) DC with Lagrange polynomials. The same Legendre polynomials of order 4 will be considered for DMS with IRK, and DC.

In the case of rigid-body dynamics, the dynamics constraint Eq.~(\ref{eq:dynamic_constraint_continuous}) of a multibody system, with generalized coordinates $\q$, velocities $\qdot$, and accelerations $\qddot$ is:
\begin{align} \label{eq:equationofmotion}
M(\q) \qddot + N(\q, \qdot) &= S^{\top} \btau_J,
\end{align}
\noindent where $M(\q)$ is the mass matrix, $N(\q,\qdot)$ is the vector of nonlinear effects and gravity, $\btau_J$ is the joint torques vector, and $S$ the selection matrix. 
By letting the state and control vectors be $\states=[\q, \qdot]$ and $\controls = \btau_J$ respectively, this equation can be rewritten as an ODE of the form~$\statesdot = \sysd(\states,\controls)$ to match the previously introduced equations.

Concerning DMS with IRK and DC, two ways to compute the dynamics constraints will be presented (Eqs.~\ref{eq:fd_constraint} and \ref{eq:id_constraint}). At each knot interval $[t_n, t_{n+1}]$, the states trajectory $\states$ is modeled by a polynomial $\bm{p}_n$, concatenating the generalized coordinates polynomial $\bm{p}_n^{\q}$ and the generalized velocities polynomial $\bm{p}_n^{\qdot}$. The generalized accelerations polynomial is directly derived from $\bm{p}_n^{\qdot}$ as $\bm{p}_n^{\qddot}(t) = \dot{\bm{p}}_n^{\qdot}(t),~\forall t$. \textit{At collocation points}, the first constraint ensures that the derivative of the generalized coordinates polynomial matches the generalized velocities polynomial:
\begin{equation}
\forall t_{n,m},~ p_n^{\qdot}(t_{n,m}, \Tilde{\states}_{n}) = \dot{p}_n^{\q}(t_{n,m}, \Tilde{\states}_{n}).
\end{equation}

The consistency of the generalized coordinates, velocity, and acceleration polynomials can be constrained through forward dynamics (FD), as written in Eq.~\ref{eq:Colloc-Conditions}, or through inverse dynamics (ID), by rearranging the terms the equations:
\begin{align}
    \forall t_{n,m},~FD\left(\bm{p}_n^{\q}(t_{n,m}, \Tilde{\states}_{n})  , \bm{p}_n^{\qdot}(t_{n,m}, \Tilde{\states}_{n})  ,\; \btau_n \right) - \bm{p}_n^{\qddot}(t_{n,m}, \Tilde{\states}_{n}) = \bm{0} \label{eq:fd_constraint},\\
    \forall t_{n,m},~ID\left(\bm{p}_n^{\q}(t_{n,m}, \Tilde{\states}_{n}), \bm{p}_n^{\qdot}(t_{n,m}, \Tilde{\states}_{n})  ,\; \bm{p}_n^{\qddot}(t_{n,m}, \Tilde{\states}_{n}) \right) - \btau_n = \bm{0} \label{eq:id_constraint}.
\end{align}
Eq.~\ref{eq:fd_constraint} amounts to matching joint accelerations with forward dynamics ($\qddot - FD(\q, \qdot, \btau) = 0$), with $FD$ the forward dynamics algorithm using the articulated body algorithm (ABA), \textit{i.e.}, without numerical inversion of the mass matrix~\citep{Featherstone1999Divide-and-conquerAlgorithm}.
Eq.~\ref{eq:id_constraint} amounts to matching torques with inverse dynamics ($\btau - ID(\q, \qdot, \qddot) = 0$), with $ID$ the inverse dynamics using the recursive Newton-Euler algorithm (RNEA, \citep{Featherstone1987RobotAlgorithms}). 
Both algorithms run in time $O(n)$, RNEA (ID) algorithm is expected to be faster since it relies on a two-pass algorithm, whereas ABA (FD) has three recursive loops.
The two formulations will be investigated in the following.

\section{Methods}

Five OCPs, involving models of increasing complexity (from 3 to 15 DoFs; Fig. \ref{fig:kinograms} and Table \ref{tab:OCPs}), were used to compare the performances of the transcription methods.
The number of \knotpoints~(Table~\ref{tab:OCPs}) was adjusted to ensure good convergence properties for each OCP, according to our experience \citep{Bailly2021Real-TimeBiomechanics, Puchaud2023OptimalityDynamics} at least 100 \knotpoints~per~second.

 As frequently done in biomechanics, the rigid-body models were either actuated by muscles or by joint torques.
The initial and final states of the OCPs were bounded and/or constrained.
All the tasks are briefly described below and detailed in \ref{appendix}. 



\begin{table*}[ht!]
\caption{Summary of the five OCPs}
\centering
\label{tab:OCPs}
\begin{adjustbox}{max width=\textwidth}
\begin{threeparttable}
\begin{tabular}{ccc | cc ccc }
\toprule
\multicolumn{3}{l}{\textbf{OCP}} & \textbf{1} &\textbf{2} & \textbf{3} & \textbf{4} & \textbf{5} \\
\midrule
\multicolumn{3}{l}{\textbf{Model}} & Leg & Arm & Upper limb & Planar human & Acrobat \\
\multicolumn{3}{l}{\# \textbf{DoFs}} & 3 & 6 & 10 & 12 & 15 \\
\multicolumn{3}{l}{\textbf{Dynamics}} & Torque-driven & Torque-driven & \makecell{ Muscle-driven \\ with residual torques} & \makecell{Torque-driven \\ with contact} & Torque-driven \\
\multicolumn{3}{l}{\# \textbf{Actuators}} & 3 & 6 & 26 + 10 & 9 & 9 \\
\multicolumn{3}{l}{\textbf{Forward Dynamics (ABA)}\textsuperscript{a} ($\mu$s)} &
$8.1$ & $9.7$ & $13$ & $12$ & $15$ \\
\multicolumn{3}{l}{\textbf{Inverse Dynamics (RNEA)}\textsuperscript{a} ($\mu$s)} &
$6.4$ & $7.5$ & $9.8$ & $ 9.3$ & $11$ \\
\multicolumn{3}{l}{\textbf{Knot points}} & 20 & 50 & 100 & 30 & 125 \\
\multicolumn{3}{l}{\textbf{Final time} (s)} & 0.25 & 0.25 & 1.0 & 0.3 & 1.5 $\pm$ 0.05\\
\bottomrule
\end{tabular}
\begin{tablenotes}
\item Notes: \textsuperscript{a} Reported dynamics computation time is averaged over 10\textsuperscript{5} function calls.
\end{tablenotes}
\end{threeparttable}
\end{adjustbox}
\end{table*}

\subsection{Rigid-body models and task}
The following models were inferred from real systems or adapted from the literature:

\begin{itemize}
  \item [\textit{OCP1}] A 3-DoF hexapod leg Robotis BIOLOID was modeled as in \citet{Trivun2017ResilientRobot}, using three hinge joints (Z-Y-Y).
  The task was to move the end-effector between two set positions while minimizing joint torques (primarily) and joint velocities (secondary).
  
  \item [\textit{OCP2}] A 4-body, 6-DoF arm model was adapted from \citet{Docquier2019ComparisonProblem}.
  The first body is connected to the base with two hinges. 
  The second and third bodies are connected through a single hinge joint. The last body is an inverted pendulum connected through a universal joint on the third body.
  A reaching task was simulated, maintaining the pendulum vertical at the beginning, and the end of the task similarly to \citet{Docquier2019ComparisonProblem}, by mainly minimizing joint torques and their time derivatives.

   \item [\textit{OCP3}] A 5-body, 10-DoF upper-limb musculoskeletal model was derived from \citet{Wu2016}. 
  The sternoclavicular joint was modeled as a universal joint. 
  The acromioclavicular and glenohumeral joints were modeled as 3-DoF ball-and-socket joints. 
  The forearm and wrist were modeled as body segments connected to the elbow via a universal joint. 
  Twenty-six muscle units represented the major muscle groups of the shoulder complex. 
  A reaching task was simulated to move the hand between two set positions, while minimizing joint torques, joint velocities, and their respective rates of change.
  
  \item [\textit{OCP4}] A 12-DoFs planar human model (83~kg, 1.67~m) was implemented in line with \citet{Serrancoli2019AnalysisProblems}. 
  The torso segment was the free-floating base with three~DoFs.
  The neck, the left, and right hips, knees, shoulders, and elbows were revolute actuated joints.
  Half a gait cycle was simulated. Joint torques and joint accelerations were minimized.
  
    \item [\textit{OCP5}] A 15-DoFs full human body was modeled (male, 18~years~old, 73.6~kg, 170.8~cm) \citep{Puchaud2023OptimalityDynamics}: 6 unactuated DoFs for the floating base, 3 for the thoracolumbar joint (flexion, lateral flexion, and axial rotation), 2 for each shoulder (plane of elevation and elevation), and 2 for the right and left hips (flexion and abduction).
    The inertial parameters were based on 95 anthropometric measures on an elite trampolinist, in line with Yeadon's anthropometric model \citep{Yeadon1990TheBody}. 
    A $1\frac{3}{4}$ somersault with one twist was simulated from a trampoline bed until a pre-landing position, similarly to \citet{Puchaud2023OptimalityDynamics}. A multi-objective function included the joint torques and the rate of change of joint velocities.

\end{itemize}

\subsection{Transcriptions of the optimal control problems} 
Five OCP transcriptions were compared, three of them belonging to direct multiple shooting schemes: 
\begin{itemize}
    \item explicit 4th-order Runge-Kutta with five intermediate steps (\ERK),
    \item implicit Runge-Kutta with rigid-body dynamics constraints based on forward dynamics (\IRKfd), 
    \item implicit Runge-Kutta with rigid-body dynamics constraints based on inverse dynamics (\IRKid),
\end{itemize}  
and the remaining two, belonging to direct collocation schemes:
\begin{itemize}
    \item rigid-body dynamics constraints based on forward dynamics (\DCfd),
    \item rigid-body dynamics constraints based on inverse dynamics (\DCid).
\end{itemize}

A fourth-order polynomial of Legendre was used as a time grid for IRK and collocation.

The planar human \textit{OCP4} was not run with implicit dynamics as it required a specific formulation of inverse dynamics to include non-acceleration constraints at contact points (a feature not available in the Bioptim library yet).

\begin{table*}[ht!]
\caption{Size of the NLPs resulting from the transcription methods (DMS vs DC), for the five OCPs
}
\centering
\label{tab:nlpsize}
\begin{adjustbox}{max width=\textwidth}
\begin{threeparttable}
\begin{tabular}{l rr rr rr rr rr}
\toprule
\multicolumn{1}{l}{\textbf{OCP}} &
\multicolumn{2}{c}{\textbf{1}} &
\multicolumn{2}{c}{\textbf{2}} &
\multicolumn{2}{c}{\textbf{3}} &
\multicolumn{2}{c}{\textbf{4}} &
\multicolumn{2}{c}{\textbf{5}} \\
\multicolumn{1}{l}{\textbf{Model}} &
\multicolumn{2}{c}{\textbf{Leg}} &
\multicolumn{2}{c}{\textbf{Arm}} &
\multicolumn{2}{c}{\textbf{Upper limb}} &
\multicolumn{2}{c}{\textbf{Planar human}} &
\multicolumn{2}{c}{\textbf{Acrobat}} \\
\midrule
Transcription &
\textbf{DMS} & \textbf{DC} &
\textbf{DMS} & \textbf{DC} &
\textbf{DMS} & \textbf{DC} &
\textbf{DMS} & \textbf{DC} &
\textbf{DMS} & \textbf{DC}\\
\midrule
states & 126 & 606 & 612 & 3012 & 2020 & 2020 & 744 & 3624 & 3780 & 18780\\
controls & 60 & 60 & 300 & 300 & 3600 & 3600 & 360 & 360 & 1125 & 1125 \\
parameters & 0 & 0 & 0 & 0 & 0 & 0 & 0 & 0 & 1 & 1 \\
equality constraints & 132 & 612 & 621 & 3021 & 2000 & 10020 & 732 & 3612 & 4880 & 18750\\
\textit{nzt} in equality jacobian & 1111 & 4734 & 10805 & 52316 & 107800 & 343950 & 24395 & 66250 & 142182 & 382370\\
inequality constraints & 0 & 0 & 0 & 0 & 0 & 0 & 32 & 32 & 0 & 0\\
\textit{nzt} in inequality jacobian & 0 & 0 & 0 & 0 & 0 & 0 & 1008 & 1008 & 0 & 0\\
Lagrangian Hessian & 762 & 1703 & 7833 & 38004 & 144318 & 364098 & 48915 & 94035 & 118149 & 291048\\
\bottomrule
\end{tabular}
\begin{tablenotes}
\item Notes: nzt stands for non-zero terms. Lagrangian Hessian refers to the cost function.
\end{tablenotes}
\end{threeparttable}
\end{adjustbox}
\end{table*}

\subsection{Multi-start approach}
All OCPs were implemented in Bioptim \citep{Michaud2022BioptimBiomechanics}, an optimal control toolbox for biomechanics that implements both DMS and DC transcriptions. 
The code is available online \citep{Puchaud2023Ipuch/dms-vs-dc:0.2.0}.
The OCPs were run on an AMD Ryzen 9 5950X processor with 128 Go RAM; OCPs were parallelized on eight threads, four OCPs running in parallel. 
The problem was solved using IPOPT with exact Hessian~\citep{Wachter2006OnProgramming}, ma57 linear solver~\citep{Duff2004MA57Systems}, and CasAdi for algorithmic differentiation~\citep{Andersson2019CasADi:Control}.
Convergence and constraint tolerances were set to 10$^{-6}$ and 10$^{-4}$, respectively.
Since solving such NLPs may result in finding local minima, each NLP was solved 100 times using random noise applied to the initial guesses, whose magnitudes were specific for each OCP (Table \ref{tab:multistart}).
All solutions were considered for further analysis.

\subsection{Analysis} \label{method:computational}

For each OCP, we compared the convergence rates, the convergence times, and the values of the cost function at convergence resulting from each transcription. We reported their medians as they give a better sense of non-uniformly distributed data.
Dynamic consistency was evaluated through the root mean square error between optimal and re-integrated generalized coordinates at final time $t_f$. 
Re-integrated solutions were obtained by solving an initial value problem, starting from the state at $t_0$, and integrating forward in time the optimized control trajectory using optimized control with the variable-step explicit Runge-Kutta method of order 8(5,3)~(\citealt{Hairer1996SolvingII}, DOP853), implemented in \textit{scipy}.
Relative and absolute tolerances were set to $1e^{-3}$ and $1e^{-6}$, respectively. This choice of tolerances was made to ensure a high level of accuracy in the re-integration.

\section{Results}

\begin{figure*}[ht!]
\includegraphics[width=1\linewidth,trim={0 5 0 0}, clip]{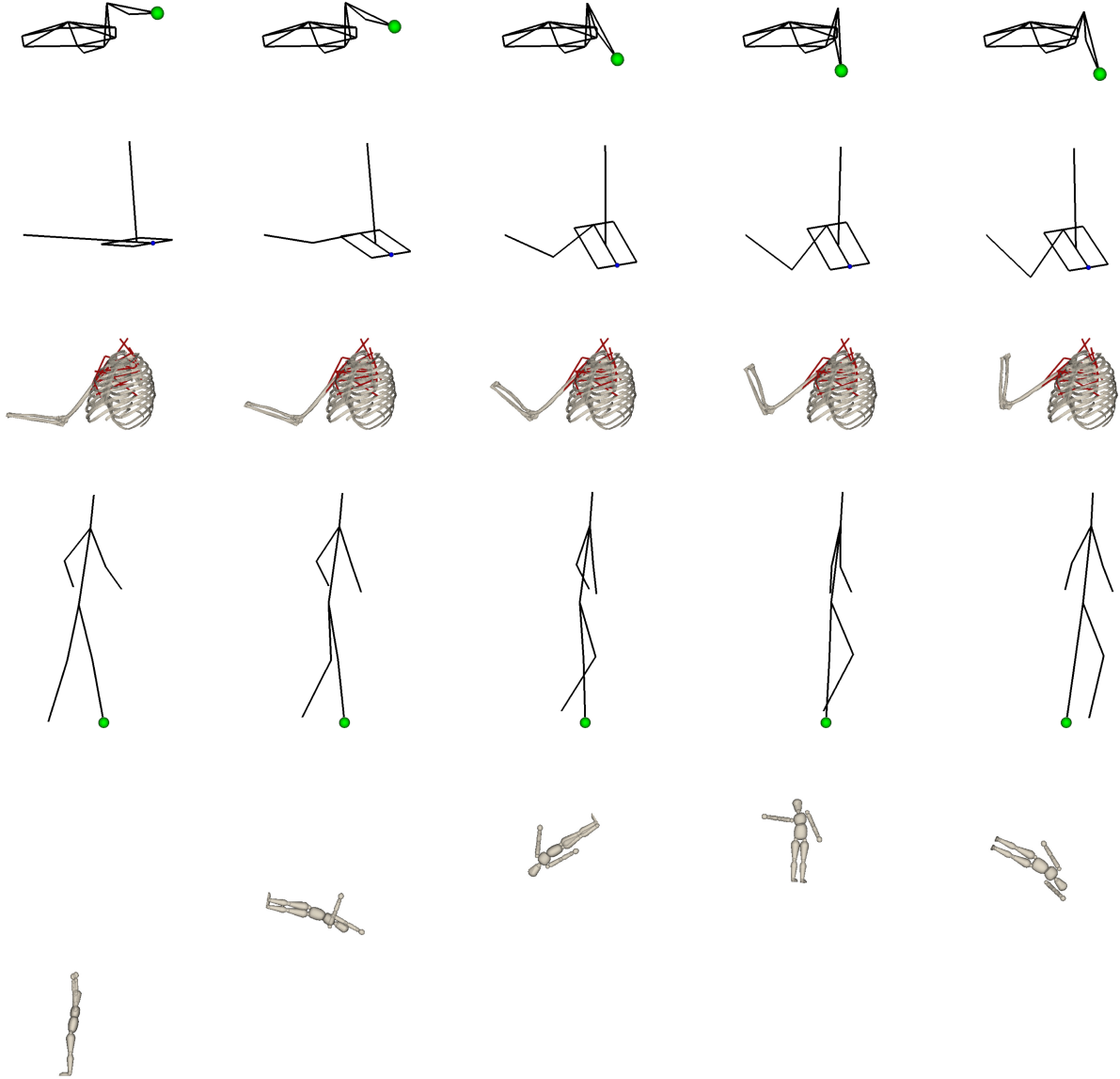}
\centering
\caption{Illustrations of results of five OCPs, from 1 to 5 (top to bottom).}
\label{fig:kinograms}
\end{figure*}

\textbf{Convergence rate.} Convergence rates for OCP1 to OCP5 were 98.2\%, 99.4\%, 100\%, 100\% and 100\% .
For OCP1, only nine simulations failed to converge (\ERK: 1, \IRKfd: 2, \IRKid: 2, \DCfd: 2, \DCid: 2) and 
only three simulations for OCP2 (\ERK: 1, \IRKfd: 1, \IRKid: 1). 
The trials that failed to converge were excluded from the subsequent analyses. 

\textbf{CPU Time.} Overall, none of the five transcriptions was faster for all OCPs. 
On the one hand, DMS transcriptions were faster than DC for OCP2 and OCP5, \IRKid~being the fastest for OCP2. 
On the other hand, DC was faster than DMS for OCP1, OCP4, and OCP3.
Among each problem, for the inverse dynamics rigid-body constraints, the median CPU times of \IRKid~and \DCid~were similar to or lower than their corresponding forward dynamics versions (\IRKfd~and \DCfd). 
For OCP1 and OCP2, \IRKid~was faster than \IRKfd~by 1.15 and 1.82 times, respectively. 
For OCP3, the two methods performed similarly. 
Also, \DCid~was faster than \DCfd~for OCP1, OCP3, and OCP5 by 1.03, 1.12, and 1.56 times respectively, except for OCP2, where \DCfd~was 1.14 times faster.

\begin{figure*}[!ht]
\includegraphics[width=0.55\linewidth,trim={0 0 0 0}, clip]{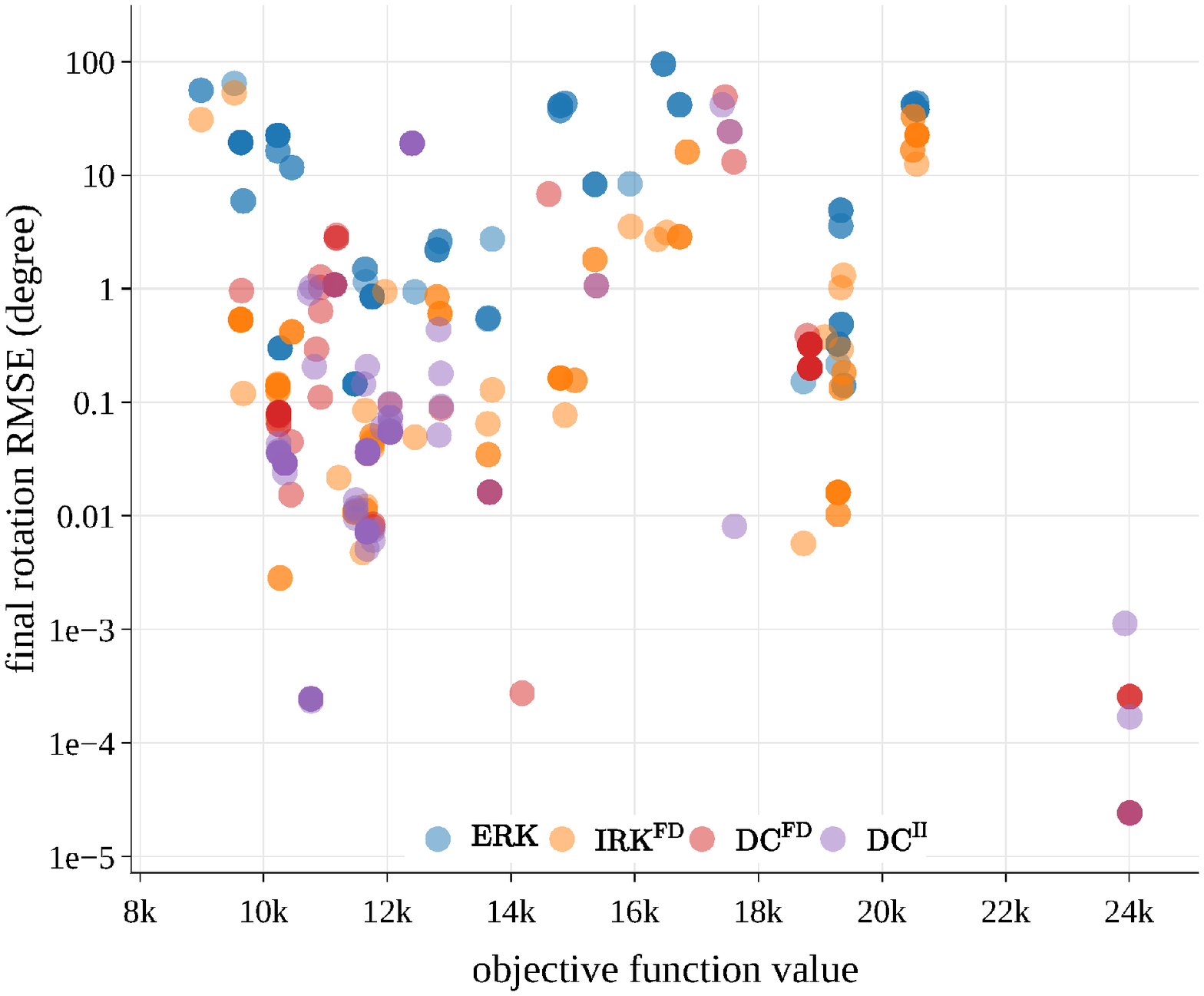}
\centering
\caption{Cost versus dynamic consistency for OCP5. The dynamics consistency y-axis has a logarithmic scale. \htmlversion}
\label{fig:costvsconsistency}
\end{figure*}

\textbf{Optimality.} All transcriptions led to similar optima.
For OCP1, two solutions in \DCfd converged to a higher local minimum (7.4 e-4 instead of 7.2e-4), and for OCP4, four solutions (2 in \ERK, 2 in \IRKfd) converged to a higher local minimum (3691 instead of 884).
The OCP2 presented two clusters of optimal solutions.
For each transcription, the first cluster contained 63\% to 68\% of the optimal solutions (\DCid and \DCfd, respectively). 
In other words, all transcription methods performed similarly.
For OCP5, the cost values were more scattered. 
\ERK, \IRKid, and \IRKfd led to similar local minima, but \DCfd showed a larger variability with a lower median.
Furthermore, \DCid led to a lower variability. 
Finally, the lowest minima were found by DMS transcriptions, but these solutions presented one of the lowest dynamic consistency compared to other local minima, see~Fig.~\ref{fig:costvsconsistency}. 

\textbf{Dynamic Consistency.} In general, \IRKfd, \IRKid, \DCfd, and \DCid~led to lower errors than \ERK~after re-integrating the optimal solution using DOP853. Errors were the highest for \ERK in OCP2, 4, and 5. 
Moreover, \ERK~led to RMSE errors up to 100\degree~by the end of the simulation of OCP5, even though these solutions were related to the lowest local minima found with all transcriptions. 
Thus, \IRKfd, \IRKid, \DCfd, and \DCid  solutions are more dynamically consistent than \ERK.

\begin{figure*}
\includegraphics[width=1\linewidth,trim={0 25 60 20},clip]{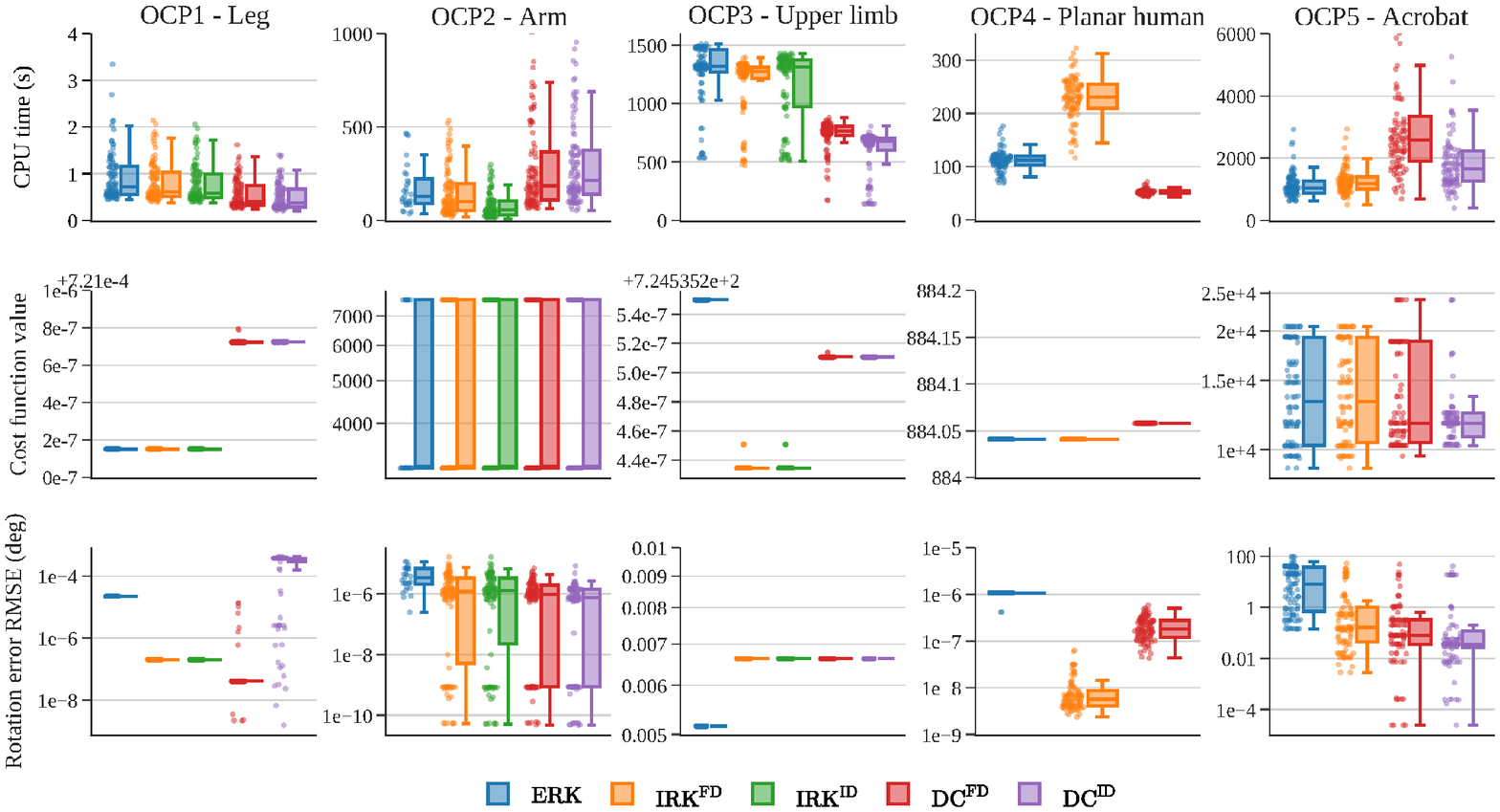}
\centering
\caption{CPU time, cost function values, and  dynamic consistency for each OCP and each transcription. The y-axes have a logarithmic scale for cost function values and dynamic consistency. \htmlversion}
\label{fig:time_iter}
\end{figure*}

\section{Discussion}
This study aimed to evaluate the performance of various direct optimal control transcriptions and dynamic formulations on five OCPs, related to rigid-body dynamics and biomechanics.
We found that the five OCPs of increasing complexity (including torque-driven, and muscle-driven dynamics) were solved in a similar amount of CPU time, whatever the transcription.
Among the transcriptions, those incorporating implicit dynamics, such as DMS with IRKs or DC, preserved dynamic consistency better. 
Rigid-body dynamics constraints based on inverse dynamics slightly hastened the convergence.

\subsection{Limitations and strengths}
We will discuss these findings in line with two main limitations.  
First, a good initial guess helps the NLP solver to converge toward an optimal solution.  
However, all examples were treated as predictive simulations without \textit{a priori} knowledge of the solution. 
Interior-point algorithms are more robust than SQP solvers to such poor initial guesses \citep{Betts2010PracticalProgramming}. 
Nevertheless, SQP solvers are well designed to solve DMS and DC quickly, Acados \citep{Verschueren2022AcadosaControl} being a good example for real-time applications. 
In the present study, we only used one interior-point solver, IPOPT, commonly used in biomechanics \citep{Dembia2019, Lee2016GeneratingMATLAB, Bailly2021OptimalModel}. 
Since performances depend on the NLP solver, our findings cannot be extrapolated to other interior-point or SQP algorithms.  

Secondly, our DMS implementation did not correspond to the state-of-the-art algorithm of \cite{Bock1984AProblems}. The DMS transcription in Bioptim is mainly based on various fixed-step ERK and IRK. In contrast, a strong feature of DMS is to rely on an adaptive step size ODE solver to control the accuracy of the dynamics without changing the size of the NLP problem.
On the contrary, DC requires mesh refinements to improve the dynamic consistency, changing the problem size.
Fixed-step ODE solvers are more convenient and efficient for algorithmic differentiation. While Bioptim (through CasADI, \citealt{Andersson2019CasADi:Control}) is connected to CVODES \citep{Serban2008CVODES:SUNDIALS}, a solver with sensitivity analysis capabilities, our experience with this solver resulted in a significant increase in CPU times.
Nevertheless, the 4th order polynomial and fixed-step size result in errors less than 0.01\degree~in 4 out of 5 OCPs, by the end of the simulation compared to the explicit Runge-Kutta method of order 8(5,3) DOP853. We believe these dynamical errors below 1\degree~are accurate enough for predictive simulation in biomechanics.
As shown in Sec.~\ref{sec:Background}, having fixed-step ODE solvers (5-step RK4 and 4th order Legendre collocation) reduces the difference between DMS and DC implementations on two aspects: (\textit{i}) where the dynamics constraints are solved (ODE solver~\textit{versus}~NLP solver, respectively); (\textit{ii}) the number of NLP variables and constraints (collocation points being bounded in DC). Consequently, this led to similar results in the present study.

Despite these two limitations, the strengths of the present study should be acknowledged.
First, we included various OCPs with increasing complexity from 3 to 15 DoFs, requiring from 0.45 to 1366 seconds to be solved. We included torque-driven and muscle-driven dynamics, and NLPs had a significant problem size with equality constraints ranging from 132 to 10 020. 
Within this diversity, no clear conclusion can be drawn, showing that the best transcription cannot be anticipated in the context of trajectory optimization for multibody systems.
According to \citet{Febrer-Nafria2022PredictiveReview} review, most studies considered the performance of their optimal control transcription based on a single transcription and one model. Even if we did not reach the complexity of some OCP available in the literature, up to 31 DoFs and 92 muscles \citep{Falisse2022ModelingWalking}, our study provides a more comprehensive perspective. 

Second, the transcriptions were all coded using the same software, Bioptim. 
All transcription methods differ by only two instructions (see code repository \citep{Puchaud2023Ipuch/dms-vs-dc:0.2.0}). 
Having a unique modeler and unique software to compare transcriptions greatly reduces the bias and possibility of errors from hand-coded scripts, lacking support for versatility formulations \citep{Lee2016GeneratingMATLAB, Koelewijn2020AEnvironment, Febrer-Nafria2021PredictionModel}. 
Additionally, the Lagrangian cost functions in both DMS and DC were computed with the same Euler forward approximation, making the software differences between DMS with IRK and DC minimal. Any variations in performance can thus be attributed to the fundamental differences between the two transcription approaches.
This lack of homogeneity in implementations led to large differences in performance in inverse approach with musculoskeletal models \citep{Kim2018SimilaritiesSystem}. 
Other software used for optimal control in biomechanics, such as Moco \citep{Dembia2019} or Muscod II \citep{Leineweber2003AnApplications} only provide a DC or a DMS solver.
Only a few software enable solving OCPs with DMS and DC, such as Horizon \citep{Ruscelli2022Horizon:Systems}, CasAdos \citep{Frey2022FastCasADi}, and Bioptim being specifically tailored multibody systems including biomechanics~\citep{Michaud2022BioptimBiomechanics}. 

\subsection{CPU time to convergence}
In the literature, predictive simulations with musculoskeletal models converge between 10 min to 100 h \citep{Febrer-Nafria2022PredictiveReview}. In our case, the upper limb musculoskeletal model took up to 25 min to converge. 
Full-body acrobat model with 15 DoFs (OCP4) had medians from 17.3 to 42.9 min CPU time and maximally reached 1.6 hours to converge.
These CPU time were quite low, probably because no contact mechanics was considered. In our case, CPU times to solve OCPs were similar for DMS and DC. 
In contrast, \cite{Porsa2016DirectOpenSim} reported that direct single shooting was 249 times slower than DC. 
It may not be appropriate to compare direct single shooting with DC or DMS in terms of computational efficiency, 
so this transcription was not considered in the present paper.
Indeed, The DMS transcription has been initially developed to leverage the issue that rises from the direct single shooting transcription \citep{Bock1984AProblems}.
In our case, CPU times were similar for DMS and DC, DC was about 5 times lower than DMS in OCP4 and 4 times higher in OCP2.
None of the transcriptions consistently showed a faster convergence among these five examples, and no other clear cause was found, other than the transcriptions being similar as presented in Sec.~\ref{sec:Background}.

Nonetheless, inverse dynamics constraints usually led to faster convergence. 
In all cases but OCP3, transcriptions based on inverse dynamics constraints (\IRKid~and \DCid~) resulted in faster convergence than forward dynamics constraints (\IRKfd~and \DCfd).
It confirms the hypothesis made in Sec.~\ref{back:dynamicconstraint} and strengthens previous findings \citep{Puchaud2023OptimalityDynamics, Ferrolho2021InverseOptimization}.
Running the forward dynamics function took around 1.3 times longer than the inverse dynamics function for our multibody models (Table~\ref{tab:OCPs}).
This ratio was never reached when comparing the time to converge between transcriptions (\IRKfd~vs \IRKfd, and \DCfd~vs \DCid).
Indeed, most of the time is spent evaluating the Lagrangian Hessian and the constraint Jacobian, not necessarily calling the dynamics function.
Still, we recommend the use of inverse dynamics defects when possible.

Finally, the initial guess has a significant impact on CPU time. 
An initial guess that converges quickly with one transcription does not guarantee that another initial guess will converge as fast. Moreover, an initial guess that converges quickly, does not guarantee similar results with another transcription.
Therefore, the variation in CPU time can be attributed more to the initial solution rather than the transcription method used.
If possible, we suggest providing a minimally data-informed solution to reduce the CPU time to convergence for predictive simulation in biomechanics, such as for OCP3, 4, and 5 in the present study.
But most importantly, to fairly compare different transcriptions, we suggest employing a multi-start approach. 

\subsection{Optimality}

Three interesting behaviors emerged from the present study's five OCP examples. In OCP1, 3, and 4, a unique optimal solution was found, whatever the transcription method and the initial guess. The cost values (Fig.~\ref{fig:time_iter}), states, and controls, were nearly identical.
The cost function was simple enough almost to admit a unique solution.  
The anecdotal differences in the cost between the methods may be explained by the slight change in dynamic consistency (cost increasing when consistency decreasing) or by the Karush–Kuhn–Tucker (KKT) conditions being different by construction between DMS and DC while the exit tolerance (1e-6) was the same. 

In OCP2, two clusters of solutions were found, one being a sub-optimum. Again, since the proportion of global \textit{versus} local optima were similar between the methods (from 63\% to 68\%), our recommendation is to use a multi-start approach whatever the transcription methods or to have an \textit{a priori} knowledge of the solution. 

In OCP5, several local minima were found for all the transcription methods. 
Finding a variety of local minima is usual in optimal control and specifically in biomechanics \citep{Falisse2019, Gidley2018PerformancePedaling}. 
Indeed, the objective function may include several terms, such as metabolic cost, muscle torque, stability, and penalty regularization. The weighting of these multi-objective functions is challenging, as they are usually manually fine-tuned \citep{Brown2020, Febrer-Nafria2021PredictionModel, Felis2016SynthesisMethods}. Thus, these complex objective functions may compensate each other, and, with different initial guesses, they may converge toward different optima. Nonetheless, whatever the transcription, local minima can be found, and local minima can be similar.

The present study does not support the finding: \citet{Docquier2019ComparisonProblem} found that rigid-body dynamics constraints based on forward dynamics lead to worse local optimal solutions than the constraints based on inverse dynamics. In OCP5, \DCid~led to lower local minima, and OCPs converged to the same optima. Further studies are needed to explore the effect of rigid-body dynamics constraints on optimality.
However, it is worth noting that the most optimal solutions showed low dynamic consistency. 
Indeed, the best solutions in terms of cost (see Fig.~\ref{fig:costvsconsistency}) presented low dynamic consistency (re-integration error higher than 1\degree~and up to 100\degree~in \ERK). Thus, evaluating only the quantity of low local minima can be misleading as the optimal solution can be dynamically inconsistent since the transcriptions may exploit its weaknesses to get better cost function value.

\subsection{Dynamic consistency}
Dynamic consistency is the criterion that showed the largest differences between the transcriptions, especially \ERK~\textit{versus} the collocation-based methods (IRK and DC).
IRK based on 4th-order Legendre polynomials provided consistent dynamics for all simulations. 
Both \IRKid~and \IRKfd~based on Legendre collocation time-grid showed reliable dynamic consistency for the presented OCPs, and, \ERK~showed the largest errors, especially for OCP5. 

IRKs are efficient when wolving OCPs with SQP methods \citep{Quirynen2015LiftedControl}. For predictive simulations, implicit integrators are recommended since they better control the integration error and, consequently, the dynamic consistency. Further analyses should be made to compare the result with adaptive step size ODE solvers to characterize DMS transcription better to exclude the use of direct multiple shooting with explicit ODE solvers for predictive simulation. These results strengthen previous recommendations to use implicit dynamics 
\citep{VanDenBogert2011ImplicitControl, DeGroote2016}, but do not exclude DMS as a transcription.
Another type of integrators may also be considered to enhance dynamic consistency, namely symplectic ODE solvers such as Leapfrog integrator \citep{Ruscelli2022Horizon:Systems}, or variational integrators that rely on discrete mechanics \citep{Johnson2009ScalableCoordinates} and may be beneficial for their energy conservation properties for biomechanical motions which last longer than 1~s (e.g., diving or repetitive tasks).

\section{Conclusion}
Our study highlights that the five transcriptions, including DMS and DC, took similar CPU time to solve the five OCPs related to rigid-body dynamics. 
To effectively compare the different transcriptions, we suggest utilizing a multi-start approach and a variety of examples.
Among the transcriptions we evaluated, those incorporating implicit dynamics, such as direct multiple shooting with IRKs or direct collocation, were particularly effective in finding optimal solutions and preserving dynamic consistency.
Additionally, we recommend employing rigid-body dynamics constraints based on inverse dynamics when available to gain extra time for convergence. We recommend using frameworks to facilitate this process and identify the most suitable transcription method for a given OCP.

\section{Acknowledgment}
The NSERC Discovery grant of Micka\"el Begon (RGPIN-2019-04978) funded this study.

\bibliographystyle{unsrtnat}
\bibliography{Mendeley.bib}  

\appendix
\newpage
\setcounter{figure}{0}
\setcounter{table}{0}

\section{Detailed OCPs} \label{appendix}
\subsection{Torque-driven robot leg pointing task}
\textbf{Multibody model.} A three-DoF kinematic chain articulated with three hinges (Fig. \ref{fig:kinograms}).

\textbf{Task.} The end-effector of the robot is constrained to reach specified locations at $t_0$ and $t_f$ with null joint velocities $\qdot(t_0)=\qdot(t_f)=\bm{0}$. A Lagrange term minimizes joint torques and velocity, leading to the following cost function:

\begin{equation} \label{eq:cost_leg}
\begin{aligned}
\Phi = \omega_1 \int_{t_0}^{t_f} \btau^2 \: dt + \omega_2 \int_{t_0}^{t_f} \qdot^2  \: dt
\end{aligned}
\end{equation}
with $\omega_1$ and $\omega_2$ set to 1 and 1e-6, respectively.

\subsection{Torque-driven 6 DoF pendulum equilibrium}
\textbf{Multibody model.}
The model is a 3D robotic arm with an inverted pendulum attached to its end-effector, adapted from \citet{Docquier2019ComparisonProblem}. 
It consists of 4 bodies and 6 DoFs. 
The first element is connected to the base with two hinges. 
A hinge joint connects the second and third elements. 
A universal joint characterizes the pendulum's motion relative to the end-effector.

\textbf{Task.}
The robot's end-effector is constrained to reach a specified location at $t_0$ and $t_f$. 
The velocity of the end-effector is also constrained to be zero.
In a stable equilibrium, the initial conditions of the system are given by: 
$\qdot(t_0)=\bm{0}$, and $\qddot(t_0)=\bm{0}$.
A terminal constraint ensures the pendulum remains vertical at the final instant $t_f$. 
The following objective is defined:

\begin{equation} \label{eq:cost_leg_annex}
\begin{aligned}
\Phi = \omega_1 \int_{t_0}^{t_f} \btau^2 \: dt 
+ \omega_2 \int_{t_0}^{t_f} \Delta\btau^2 \: dt 
+ \omega_3 \int_{t_0}^{t_f} \qdot^2  \: dt
\end{aligned}
\end{equation}
with $\Delta\btau$, the rate of change of the torques $\btau$. 
Weights $\omega_1$, $\omega_2$, $\omega_3$, are set to 1, 1 and 1e-2, respectively.

\subsection{Muscle-driven upper limb pointing task}
\textbf{Multibody model.} A 5-segment, 10-DoF upper-limb musculoskeletal model from \citep{Wu2016} is implemented in Biorbd. 
The forearm and wrist are modeled as rigid-body segments connected to the elbow via a 2-DoF universal joint, and the glenohumeral joint and acromioclavicular joints are modeled as 3-DoF ball-and-socket joints.
The sternoclavicular joint is modeled as a 2-DoF joint. 
Twenty-six muscle units represented the major muscle groups of the axial scapula, axial humerus, and glenohumeral. 
Residual torques are considered for each DoF except for the glenohumeral joint. 

\textbf{Task.} The goal is to reach the head from a reference position with zero velocity at the beginning and the end of the task. 
Reference trajectories $\q$ are computed using inverse kinematics on experimental data. 
The reference trajectories bound the initial and final time of the simulation. 
Joint torques, their rate of change, and muscle forces are minimized. 
Regularization objectives are set to ease the convergence. 
Generalized coordinates of the sternoclavicular and acromioclavicular joints are minimized to maintain the scapula as close as possible to its original position. 
Joint velocities and their rate of change are also minimized. 
Finally, a Mayer term minimizes the velocity to reach zero velocity at the end. 
As a result, the cost function was a sum of Lagrange and Mayer terms:

\begin{align}
\begin{split}
    \mathcal{J} = & \:
    \omega_1 \int_{t_0}^{t_f} \btau^\top_{res} \btau_{res}\: dt 
    + \omega_2 \int_{t_0}^{t_f} \Delta\btau^\top_{res} \Delta\btau_{res} \: dt\\
    &+ \omega_3 \sum_{k=0}^4 \int_{t_0}^{t_f} q^\top_i q_i \: dt + \omega_4 \int_{t_0}^{t_f} \qdot^\top \qdot \: dt + \omega_5 \int_{t_0}^{t_f} \Delta\qdot^\top \Delta\qdot \: dt \\
    & \:  \omega_6 \: \qdot(t_f)^2,
\end{split}\\
\end{align}
where weights $\omega_1$, $\omega_2$, $\omega_3$, $\omega_4$, $\omega_5$, $\omega_6$ are respectively set to
10, 1500, 1000, 200, 150, and 0.5.

\subsection{Torque-driven Planar human}
\textbf{Multibody model.}
A 12-DoFs planar human model is implemented as in~\citet{Serrancoli2019AnalysisProblems}. 
The model has three unactuated DoFs at the torso segment to freely move in the global frame. 
It comprises two hips, two knees, a neck, two shoulders, and two elbow joints. 
The model is based on a human of 83 kg and 1.67 m in height.

\textbf{Task.}
The optimal control problem simulates a half-cycle walking gait with constrained dynamics \citep{Felis2016SynthesisMethods}.
A constraint is set to ensure the location of the right foot at $t_0$ and its zero velocity.
A constraint that defines a step length of 0.8 m enforces the initial and final location of the left foot.
A heel clearance constraint is added to ensure the lifting of the foot between 4 and 6 cm over the floor in the middle of the simulation \citep{Mariani2012HeelSensors}.
Arm swing is also enforced through initial and terminal bounds. 
The vertical component of the ground reaction force is constrained to be positive. 
Torque actuators are minimized to accomplish the task. A stability criterion minimizes the floating base and head joint accelerations with a weight of 0.01 as in previous studies~\citep{Dorn2015PredictiveWalking, Ong2019PredictingSimulations, Nguyen2019BilevelGait, Veerkamp2021EvaluatingGait}. 
Two Mayer cost functions at $t_0$ and $t_f$ ensure the constant gait velocity, by tracking the horizontal center of mass velocity at $\dot{c}_y^{ref} = 1.3 \; m/s$.
As a result, the cost function is a sum of Lagrange and Mayer terms:

\begin{align}
\begin{split}
    \Phi = & \:
    \omega_1 \int_{t_0}^{t_f} \btau^\top_J \btau_J \: dt 
    + \omega_2 \sum_{k=0}^3 \int_{t_0}^{t_f} \ddot{q}^\top_k \ddot{q}_k \: dt \\ 
    +  & \: \omega_3 \left( \left( \dot{c}_y(t_0) - \dot{c}_y^{ref} \right)^2 + \left( \dot{c}_y(t_f) - \dot{c}_y^{ref} \right)^2 \right)
\end{split}
\end{align}
where $k^{\textrm{th}} \: q$ stands for the four first DoFs. $t_0$ is the initial instant of time of the OCP. 
Weight $\omega_1$, $\omega_2$, $\omega_3$, is respectively set to
1, 0.01, and 1000. 


\textbf{Equations of motions.} The equations of motion includes an inelastic contact constraint such that:
\begin{align}
M(\q) \qddot + N(\q, \qdot) &= S^{\top} \btau_J
+ K(\q)^\top \boldsymbol{\lambda} \label{eq:equationofmotion_constrained}\\ 
K(\q) \qddot + \dot{K}(\q) \qdot &= \boldsymbol{0} \label{eq:equationofmotion_constraint}
\end{align}
where $K(\q)$ is the jacobian constraint, $\boldsymbol{\lambda}$ are the Lagrange multipliers, which correspond to the contact forces.
The system can be rewritten to be handled as an ordinary differential equation by writing Eq.~(\ref{eq:equationofmotion_constrained})~and~(\ref{eq:equationofmotion_constraint}) into a linear system of the unknowns~$\qddot$,~and~$\boldsymbol{\lambda}$:

\begin{equation}
    \begin{bmatrix} M(\q) & K(\q)^\top\\ K(\q) & \boldsymbol{0} \end{bmatrix}
    \begin{pmatrix} \qddot \\ \boldsymbol{\lambda} \end{pmatrix} =
    \begin{pmatrix} - N(\q,\qdot) + S^{\top} \btau_J \\ -\dot{K}(\q) \qdot \end{pmatrix}
    \label{eq:ID_full}     
\end{equation}

Consistency conditions ensure the first frame to respect the constraints $h(\q)=\boldsymbol{0}$ and $K(\q) \qdot = \boldsymbol{0}$. 
This equation can be re-written as a differential algebraic equation and can be integrated an ODE~$\statesdot = f(\states,\controls)$ as for the classical forward dynamics equation by setting states $\states=[\q, \qdot]$ and controls $\controls = \btau_J$.

\textbf{Bounds and initial Guess.}
The state's boundaries are set to get plausible joint ranges of motion for gait.
The torso rotational DoF could only bend forward. 
The initial and final head joint velocity is constrained at a zero velocity. The initial knee flexions are limited to $\frac{\pi}{8}$. 
To get shoulder oscillations, initial right and left shoulder states are set to $-\frac{\pi}{6}$ and $\frac{\pi}{6}$ and reversed for the final states. 
The initial and final shoulder velocities are set to zero. 
All torque controls were bounded between -500 and 500 Nm except for the three controls of the torso segment, considered the root segment, which were bounded to zero.

Initial and final poses are computed through inverse kinematics, ensuring right and left feet are at the specified locations, and root segment has positive generalized coordinates to be within boundaries. 
The initial guess of the generalized coordinates is computed from a linear interpolation from the initial to the final pose. 
The initial guesses of generalized velocities $\Dot{q}$ and control torques $\tau$ is set to zero.

\subsection{Torque-driven double twisted somersault}
\textbf{Multibody model.}
The human body model includes 15 DoFs \citep{Puchaud2023OptimalityDynamics}: 6 DoFs for the floating base, 3 for the thoracolumbar joint (flexion, lateral flexion, and axial rotation), 2 for each shoulder (plane of elevation and elevation), and 2 modeling right and left hips (flexion and abduction). 
The model inertial parameters is based on 95 anthropometric measures of an elite trampolinist (male, 18 years old, 73.6 kg, 170.8 cm), in line with Yeadon's anthropometric model \citep{Yeadon1990TheBody}.

\textbf{Task.}
The goal is to execute 1 twist in~$1\frac{3}{4}$ of somersaults from takeoff to the prelanding phase.
Final time $t_f$ is an optimized parameter of the OCP, constrained to remain between $\pm$ 0.05 \si{\second} representative of the performance of an athlete measured in~\cite{Venne2022}.
The objective function enforces the execution of the twisting somersault while minimizing hands and feet trajectories in pelvis frames. 
Integral of joint torques $\btau_j$ are minimized. 
To ensure smooth joint kinematics, joint torques, and joint velocities rate of change denoted, $\Delta \btau_J$, $\Delta \qdot_J$ are minimized.
An objective term includes states of core DoFs, namely thoracolumbar joint and hips, which were denoted $\mathcal{C}_{dof}$ so that the gymnast stays as straight as possible to comply with the sports regulations.
Specific Mayer term $\mathcal{M}$ is added to minimize angular momentum $\am$ at the initial instant of time $t_0$.
Finally, the objective function is a sum of terms $\mathcal{J}$ and $\mathcal{M}$:

\begin{align}
\begin{split}
    \mathcal{J} = & \:
    \omega_1 \sum_{k=1}^{3} \left(\Delta \underline{P}^{k}\right)^2
    + \omega_2 \int_{t_0}^{t_f} \btau^\top_J \btau_J\: dt \\
    & + \omega_2 \int_{t_0}^{t_f} \Delta \btau^\top_J \Delta \btau_J\: dt 
    + \omega_3 \int_{t_0}^{t_f} \Delta \qdot^\top_J \Delta \qdot_J\: dt \\
    & + \omega_4 \sum_{k \in \mathcal{C}_{dof}} \int_{t_0}^{t_f} {q_k}^2  \: dt +
    \omega_5 \; t_f
\end{split}\\
    \mathcal{M} = &  \:\omega_6 \; \sigma_x(t_0) + \omega_6 \; ( \sigma_y(t_0)^2 + \sigma_z(t_0)2)
\end{align}
where $\Delta \underline{P}^k$ is $k^{\textrm{th}}$ rate of change of path of hands ($k= 1,2$) or foot ($k= 3$) in pelvis frame for each instant of time of the problem, $t_0$ is the initial instant of time of the OCP. 
Weight $\omega_1$, $\omega_2$, $\omega_3$, $\omega_4$, $\omega_5$, $\omega_6$ are set to
10, 100, 100, 10, 10$^{-6}$ and 50, respectively.\\

\newpage
\section{Noise applied to initial guesses of OCP}
\begin{table*}[ht!]
\caption{Noise magnitude applied on decision variables for each OCP in percentage (\%) of bound magnitude}
\centering  
\label{tab:multistart}
\begin{threeparttable}
\begin{tabular}{c| cc  cc  ccc }     
\toprule
\textbf{OCP} & \textbf{1} &\textbf{2} & \textbf{3}& \textbf{4} & \textbf{5} \\
\midrule
\textbf{Model} & Leg & Arm & Acrobat & Planar human & Upper limb \\
\midrule
$\q_{\text{init}}$ & zeros & zeros & custom & linear & linear\\
$\qdot_{\text{init}}$ & zeros & zeros & custom & zeros & zeros \\
$\btau_{\text{init}}$ & zeros & zeros & custom & zeros & zeros\\
$\bm{a}_{\text{init}}$ & - & - & - & - & constant\\
\midrule
\multicolumn{6}{c}{Noise magnitude (\%)} \\
\midrule
$\q$ & 100 & 50 & 20 - 56 & 10 & 1 \\
$\qdot$ & 100 & 50 & 2 & 1 & 0.2 \\
$\btau$ & 20 & 50 & 2 & 1 & 1 \\
$\bm{a}$ & - & - & - & - & 2 \\
\bottomrule                                                     
\end{tabular}
    \begin{tablenotes}
      \item Notes: All the noise magnitudes were chosen to keep the initial guess within the bounds of the decision variables. 
  \end{tablenotes}
\end{threeparttable}
\end{table*}

\end{document}